\documentclass[runningheads]{llncs}

\usepackage[usenames,dvipsnames,table]{xcolor}
\usepackage{array, multirow, enumerate}
\usepackage{dsfont}  \let \mathbb=\mathds
\usepackage{amsmath, amssymb, amsfonts, amsmath, bm, stmaryrd}
\usepackage{graphicx, psfrag, tikz}

\renewcommand{\hat}{\widehat}
\newcommand{\X}{\mathbb{X}}

\newcommand{\B}{\mathbb{B}}
\newcommand{\R}{\mathbb{R}}
\renewcommand{\P}{\mathbb{P}}
\newcommand{\XX}{\mathcal{X}}
\newcommand{\YY}{\mathcal{Y}}

\newcommand{\E}{\mathbb{E}}
\newcommand{\one}{\mathds{1}}

\newcommand \EI     {\rho}

\newcommand \dy {\mathrm{d}y}

\DeclareMathOperator \argmax {{\rm argmax}}
\DeclareMathOperator \argmin {{\rm argmin}}

\newcommand \Xv  {\underline{X\mskip -1.5mu} \mskip 1mu}

\newcommand \opt  {{\rm opt}}

\newcommand \yUpp {y^{\text{upp}}}

\newcommand \Hn {H_{n}}
\newcommand \Hnc {\Hn^{\rm c}}
\newlength{\vSki} 

\begin{document}
\title{User preferences in Bayesian multi-objective optimization: the expected weighted hypervolume improvement criterion}
\titlerunning{User preferences in Bayesian multi-objective optimization}
\author{Paul Feliot\inst{1} \and
Julien Bect\inst{2} \and
Emmanuel Vazquez\inst{2}}
\authorrunning{P. Feliot et al.}
\institute{Safran Aircraft Engines, \email{paul.feliot@safrangroup.com} \and
Laboratoire des Signaux et Syst{\`e}mes (L2S), Centrale-Sup{\'e}lec, \email{firstname.lastname@centralesupelec.fr}}
\maketitle 
\begin{abstract}
  In this article, we present a framework for taking into account user
  preferences in multi-objective Bayesian optimization in the case
  where the objectives are expensive-to-evaluate black-box
  functions. A novel \emph{expected improvement} criterion to be used
  within Bayesian optimization algorithms is introduced. This
  criterion, which we call the \emph{expected weighted hypervolume
    improvement} (EWHI) criterion, is a generalization of the popular
  \emph{expected hypervolume improvement} to the case where the
  hypervolume of the dominated region is defined using an absolutely
  continuous measure instead of the Lebesgue measure. The EWHI
  criterion takes the form of an integral for which no closed form
  expression exists in the general case. To deal with its computation,
  we propose an importance sampling approximation method. A sampling
  density that is optimal for the computation of the EWHI for a
  predefined set of points is crafted and a sequential Monte-Carlo
  (SMC) approach is used to obtain a sample approximately distributed
  from this density. The ability of the criterion to produce
  optimization strategies oriented by user preferences is demonstrated
  on a simple bi-objective test problem in the cases of a preference
  for one objective and of a preference for certain regions of the
  Pareto front.

\end{abstract}
\section{Introduction}

In this article, we present a Bayesian framework for taking into
account user preferences in multi-objective optimization when
evaluation results for the functions of the problem are obtained using
a computationally intensive computer program. Such a setting is
representative of engineering problems where finite elements analysis
or fluid dynamics are used. The number of runs of the computer program
that can be afforded is limited and the objective is to build a
sequence of observation points that rapidly provides a ``good''
approximation of the set of \emph{Pareto optimal solutions}, where
``good'' is measured using some user-defined loss function.

To this end, we formulate an \emph{expected improvement} (EI)
criterion (see, e.g., \cite{jones1998efficient}) to be used within the
BMOO algorithm of~\cite{feliot2017bayesian} that uses the
\emph{weighted hypervolume indicator} (WHI) introduced
by~\cite{zitzler2007hypervolume} as a loss function. This new
criterion, which we call the \emph{expected weighted hypervolume
  improvement} (EWHI) criterion, can be viewed as a generalization of
the \emph{expected hypervolume improvement} (EHVI) criterion
of~\cite{emmerich2006single} that enables practitionners to tailor
optimization strategies according to user preferences.

The article is structured as follows. First, we recall in
Section~\ref{sec:bayesian} the framework of Bayesian
optimization. Then, we detail in Section~\ref{sec:ewhi} the
construction of the EWHI criterion and discuss computational
aspects. The ability of the criterion to produce optimization
strategies according to user preferences is then demonstrated on a
simple bi-objective test problem in the cases of a preference for one
objective and of a preference for certain regions of the Pareto front
in Section~\ref{sec:experiments}. Finally, conclusions and
perspectives are drawn in Section~\ref{sec:conclusions}.

\section{Bayesian optimization}
\label{sec:bayesian}
\subsection{The Bayesian approach to optimization}

Consider a continuous optimization problem $\mathcal{P}$ defined over
a search space $\X \subset \R^d$ and let
$\Xv = (X_1,\,X_2,\, X_3\ldots)$ be a sequence of observation points
in $\X$. The problem $\mathcal{P}$ can be, for example, an
unconstrained single-objective optimization problem or a constrained
multi-objective problem. The quality at time $n>0$ of the sequence
$\Xv$ viewed as an approximate solution to the optimization problem
$\mathcal{P}$ can be measured using a positive loss function
\begin{equation}
  \label{eq:loss-fun}
  \varepsilon_n : \Xv \mapsto \R^+\,,
\end{equation}
such that $\varepsilon_n(\Xv) = 0$ if and only if the set
$(X_1, \ldots, X_n)$ solves $\mathcal{P}$ and, given two optimization
strategies $\Xv_1$ and $\Xv_2$,
$\varepsilon_n(\Xv_1) < \varepsilon_n(\Xv_2)$ if and only if $\Xv_1$
offers a better solution to $\mathcal{P}$ than $\Xv_2$ at time
$n$. Under this framework, one can formulate the notion of improvement
as a measure of the loss reduction yielded by the observation of a new
point $X_{n+1}$:
\begin{eqnarray}
  \label{eq:improvement}
  I_{n+1} &=& \varepsilon_{n}(\Xv) - \varepsilon_{n+1}(\Xv)\,,\, n \geq 0\,.
\end{eqnarray}
The improvement is positive if $X_{n+1}$ improves the quality of the
solution at time $n+1$ and zero otherwise.

Assume a statistical model with a vector-valued stochastic process
model $\xi$ with probability measure $\P_0$ representing prior
knowledge over the functions involved in the optimization problem
$\mathcal{P}$. Under the Bayesian paradigm, optimization algorithms
are crafted to achieve, on average, a small value of
$\varepsilon_n(\Xv)$ when $n$ increases; where the average is taken
with respect to $\xi$. In this framework, the choice of the
observation points $X_i$ is a sequential decision problem. The
associated Bayesian-optimal strategy for a finite budget of $N$
observations is, however, not tractable in the general case for $N$
larger than a few units. To circumvent this difficulty, a common
approach is to consider one-step look-ahead strategies (also referred
to as myopic strategies, see, e.g.,~\cite{kushner1964new,mockus78}
and~\cite{benassi2013nouvel,GinsbLeRiche2009} for discussions about
two-step look-ahead strategies) where observation points are chosen
one at a time to minimize the conditional expectation of the future
loss given past observations:
\begin{eqnarray}
  \label{eq:oneStepLookAhead}
  X_{n+1} &=& \argmin_{x\in\X} \E_n\bigl( \varepsilon_{n+1}(\Xv)
  \mid X_{n+1}=x \bigr) \nonumber\\
  		  &=& \argmax_{x\in\X} \E_n\bigl( \varepsilon_{n}(\Xv) - \varepsilon_{n+1}(\Xv) \mid X_{n+1}=x \bigr) \nonumber\\
  		  &=& \argmax_{x\in\X} \E_n\bigl( I_{n+1}(\Xv) \mid X_{n+1}=x \bigr)\,,\, n \geq 0\,,
\end{eqnarray}
where $\E_{n}$ stands for the conditional expectation with respect to
$X_1,\, \xi(X_1),\,\ldots$, $X_n,\,\xi(X_n)$. The function
\begin{equation}
  \label{eq:EI_mono}
  \EI_n : x \mapsto \E_n \bigl( I_{n+1}(\Xv) \mid X_{n+1}=x \bigr) \,,\, n \geq 0\,,
\end{equation}
is called the \emph{expected improvement}~(EI). It is a popular
sampling criterion in the Bayesian optimization literature for
designing optimization algorithms (see,
e.g.,~\cite{jones1998efficient,schonlau1998global} for applications to
constrained and unconstrained global optimization problems).

\subsection{Multi-objective Bayesian optimization} 
We focus in this work on unconstrained multi-objective optimization
problems. Given a set of objective functions $f_j:\X\to\R$,
$j = 1,\, \ldots,\, p$, to be minimized, the objective is to build an
approximation of the Pareto front and of the set of corresponding
solutions
\begin{equation} 
  \Gamma = \{x \in \X : \nexists\, x'\in\X \text{\ such that \ }
  f(x') \prec f(x) \}\,, 
\end{equation}
where $\prec$ stands for the Pareto domination rule defined on $\R^p$ by
\begin{equation}
  \label{eq:pareto-dom}
  y=(y_1,\,\ldots,\,y_p) \prec z=(z_1,\,\ldots,\,z_p) \Longleftrightarrow \left\{
    \begin{array}{l  l}
      \forall i \le p, &\; y_i \leq z_i\,, \\[0.6em]
      \exists j \le p, &\; y_j < z_j\,.
    \end{array}
  \right.  
\end{equation}

In this setting, it is common practice to measure the quality of
optimization strategies using the hypervolume loss function (see,
e.g.,
\cite{knowles2002metrics,laumanns1999approximating,zitzler1998multiobjective})
defined by
\begin{equation}
  \label{eq:hypervolume_loss}
  \varepsilon_n(\Xv) = \left| H \setminus H_n \right|\,,  
\end{equation}
where $\left| \,\cdot\, \right|$ denotes the usual (Lebesgue) volume
measure in~$\R^p$ and where, given an upper-bounded set $\B$ of the
form~$\B = \left\{ y \in \R^p;\; y \le \yUpp \right\}$ for
some~$\yUpp \in \R^p$, the subsets
\begin{equation}
H = \{y \in \B\,;\, \exists x \in \X\,,\, f(x) \prec y\}\,,  
\end{equation}
and
\begin{equation}
H_n = \{y \in \B\,;\, \exists i\le n\,,\, f(X_i) \prec y\}\,,  
\end{equation} 
denote respectively the subset of points of $\B$ dominated by the
points of the Pareto front and the subset of points of $\B$ dominated
by~$\left(f(X_1),\ldots,f(X_n)\right)$. The set $\B$ is introduced to
ensure that the volumes of~$H$ and~$H_n$ are finite.

Using the loss function~\eqref{eq:hypervolume_loss}, the improvement
function~\eqref{eq:improvement} takes the form
\begin{equation}
  I_{n+1}\left( \Xv \right) 
  = \left| H \setminus H_n \right| - \left| H \setminus H_{n+1} \right|
  = \left| H_{n + 1} \setminus H_n \right|\,,
\end{equation} 
and an expected improvement criterion can be formulated as
\begin{eqnarray}
  \label{eq:EI-multi}
  \EI_{n}(x) &=& \displaystyle \E_n \bigl( I_{n+1}(\Xv) \mid X_{n+1}=x \bigr) \nonumber \\ 
  &=& \displaystyle \mathbb{E}_{n}\left(\int_{\B \setminus H_n}
    \mathds{1}_{\xi (x)\prec y}\, \dy\right) \nonumber \\
  &=& \displaystyle \int_{\B \setminus H_n}
  \mathbb{E}_{n} \left(\mathds{1}_{\xi (x)\prec y}\right)\, \dy
  \nonumber \\
  &=& \displaystyle \int_{\B \setminus H_n} \P_n\left(\xi(x)\prec y\right)\,
  \dy\,,
\end{eqnarray}
where $\P_n$ stands for the probability $\P_0$ conditioned on
$X_1, \xi(X_1), \ldots, X_n,\xi(X_n)$.  The multi-objective sampling
criterion~\eqref{eq:EI-multi} is called the \emph{expected hypervolume
  improvement}~(EHVI) criterion. It has been proposed and studied by
Emmerich and
coworkers~\cite{emmerich2005,emmerich2006single,emmerich2008computation}.

\section{Expected weighted hypervolume improvement (EWHI)}
\label{sec:ewhi}
\subsection{Formulation of the criterion}

To measure the quality of Pareto approximation sets according to user
preferences, Zitzler et al. (2007) proposed to use a user-defined
continuous measure in the definition of the hypervolume
indicator\footnote{In the original definition, the authors introduce
  additional terms to weight the axis. In this work, one of our
  objective is to get rid of the bounding set $\B$, as proposed by
  \cite{emmerich2014reference}. Therefore we do not consider these
  terms.} instead of the Lebesgue
measure~(see~\cite{zitzler2007hypervolume}):
\begin{equation}
  \label{eq:weighted_loss}
  \varepsilon_n(\Xv) = \mu(H \setminus H_n)\,,
\end{equation}
where the measure $\mu$ is defined by $\mu(\dy) = \omega(y)\,\dy$
using a positive weight function $\omega: \R^p \rightarrow \R^+$. The
value $\omega(y)$ for some $y \in \R^p$ can be seen as a reward for
dominating $y$ that the user may specify. Optimization strategies
crafted using the loss function~\eqref{eq:weighted_loss} have been
studied by~\cite{auger2009articulating,auger2009investigating,emmerich2014reference,zitzler2007hypervolume}.

Observe that, as discussed by~\cite{emmerich2014reference}, assuming
that $\mu$ possesses the bounded improper integral property,
\eqref{eq:weighted_loss} is well defined and upper-bounding values are
no longer required in the definition of the sets $H$ and $\Hn$, which
can be redefined as:
\begin{equation}
  \left\{
    \begin{array}{lcl}
      H &=& \{y \in \R^p\,;\, \exists x \in \X\,,\, f(x) \prec y\}\,,    \\[8pt]
      H_n &=& \{y \in \R^p\,;\, \exists i\le n\,,\, f(X_i) \prec y\}\,.
    \end{array}
  \right.
\end{equation}

Similarly to~\eqref{eq:hypervolume_loss}, the improvement function
associated to the loss function~\eqref{eq:weighted_loss} takes the
form
\begin{equation}
  I_{n+1}\left( \Xv \right) 
  = \mu(H \setminus H_n) - \mu(H \setminus H_{n+1})
  = \mu(H_{n+1}\setminus H_{n})\,,
\end{equation} 
and an expected improvement criterion can be formulated as:
\begin{eqnarray}
  \label{eq:EWHI}
  \EI_{n}(x) &=& \displaystyle \E_n \bigl( I_{n+1}(\Xv) \mid X_{n+1}=x \bigr) \nonumber \\ 
  &=& \displaystyle \mathbb{E}_{n}\left(\int_{\Hnc}
    \mathds{1}_{\xi (x)\prec y}\, \mu(\dy)\right) \nonumber \\
  &=& \displaystyle \int_{\Hnc} \P_n\left(\xi(x)\prec y\right)\, \omega(y)\, \dy\,,
\end{eqnarray}
where $\Hnc$ denotes the complementary of $H_n$ in $\R^p$. By analogy
with the EHVI criterion, we call the expected improvement
criterion~\eqref{eq:EWHI} the \emph{expected weighted hypervolume
  improvement} (EWHI) criterion.

\subsection{Computation of the criterion}
Under the assumption that the components $\xi_i$ of $\xi$ are mutually
independent stationary Gaussian processes, which is a common modeling
assumption in the Bayesian optimization literature (see, e.g.,
\cite{santner2003design}), the term $\P_n\left(\xi(x)\prec y\right)$
in the expression~\eqref{eq:EWHI} of the EWHI can be expressed in
closed form: for all $x \in \X$ and $y \in \Hnc$,
\begin{equation}
  \label{eq:closedFormIntegrand}
  \displaystyle \P_n\left(\xi(x) \prec y\right) = \displaystyle
  \prod_{i=1}^p
  \Phi \left( \frac{y_i- \hat\xi_{i,n}(x)}{\sigma_{i,n}(x)} \right)\,,
\end{equation}
where~$\Phi$ denotes the Gaussian cumulative distribution function
and~$\hat\xi_{i,n}(x)$ and $\sigma^2_{i,n}(x)$ denote respectively the
kriging mean and variance at~$x$ for the~$i^\text{th}$ component
of~$\xi$ (see, e.g., \cite{santner2003design,williams2006gaussian}).

The integration of~\eqref{eq:closedFormIntegrand} over $\Hnc$ on the
other hand, is a non-trivial problem. Besides, it has to be done
several times to solve the optimization
problem~\eqref{eq:oneStepLookAhead} and choose $X_{n+1}$. To address
this issue, we propose to choose $X_{n+1}$ among a set of predefined
candidate points obtained using sequential Monte-Carlo techniques as
in~\cite{feliot2017bayesian}, and derive a method to compute
approximations of~\eqref{eq:EWHI} with arbitrary weight functions
$\omega$ for this set.

Let then
$\XX_n = \left(x_{n,k}\right)_{1 \leq k \leq m_x} \in \X^{m_x}$ be a
set of $m_x$ points where $\EI_n$ is to be evaluated and denote
\begin{equation}
  \label{eq:ei_int}
  \EI_{n,k} = \EI_{n}(x_{n,k}) =
  \int_{\Hnc} \omega(y)\, \P_n \left( \xi(x_{n,k}) \prec y \right)\, \dy \,,\quad 1 \leq k \leq m_x\,.
\end{equation}

Using a sample $\YY_n = (y_{n,i})_{1 \leq i \leq m_y}$ of $m_y$ points
obtained from a density $\pi_n$ on $\Hnc$ with un-normalized density
$\gamma_n$ and with normalizing constant
\begin{equation}
  Z_n = \int_{\Hnc} \gamma_n(y)\, \dy\,,
\end{equation}
an importance sampling approximation of the
$(\EI_{n,k})_{1 \leq k \leq m_x}$ can be written as
\begin{equation}
  \label{eq:ei_int_approx}
  \hat\EI_{n,k} = \frac{Z_n}{m_y} \sum_{i=1}^{m_y} \frac{\omega(y_{n,i})\, \P_n \left( \xi(x_{n,k}) \prec y_{n,i} \right)}{\gamma_n(y_{n,i})}\,,\quad 1 \leq k \leq m_x\,.
\end{equation}  

To obtain a good approximation for all $\hat\EI_{n,k}$ using a single
sample $\YY_n$, the un-normalized density $\gamma_n$ can be chosen to
minimize the average sum of squared approximation errors:
\begin{equation}
  \begin{array}{l}
    \displaystyle \E \left( \sum_{k=1}^{m_x} \left( \hat\EI_{n,k} - \EI_{n,k} \right)^2 \right) \\
    = \displaystyle \frac{1}{m_y} \sum_{k=1}^{m_x} \left( Z_n \int_{\Hnc} \frac{\omega(y)^2\,\P_n \left( \xi(x_{n,k}) \prec y \right)^2}{\gamma_n(y)^2} \gamma_n(y)\, \dy - \EI_{n,k}^2 \right) \\
    = \displaystyle \frac{1}{m_y} \left( Z_n\int_{\Hnc} \frac{\sum_{k=1}^{m_x} \omega(y)^2\, \P_n \left( \xi(x_{n,k}) \prec y \right)^2}{\gamma_n(y)^2} \gamma_n(y)\, \dy - \sum_{k=1}^{m_x} \EI_{n,k}^2 \right) \,.
  \end{array}
\end{equation} 

This leads, using the Cauchy-Schwarz inequality (see, e.g.,
\cite{bect2015echantillonnage}), to the definition of the following
density on $\Hnc$:
\begin{equation}
  \label{eq:l2opt}
  L_2^{\opt}(y) \propto \gamma_n(y) = \displaystyle \sqrt{\sum_{k=1}^{m_x} \omega(y)^2\, \P_n \left( \xi(x_{n,k}) \prec y \right)^2}\,.
\end{equation}

To obtain a sample distributed from the $L_2^{\opt}$ density and carry
out the approximate computation of the EWHI
using~\eqref{eq:ei_int_approx}, we resort to sequential Monte-Carlo
techniques as well~(see, e.g.,
\cite{au2001estimation,del2006sequential,feliot2017bayesian}). The
algorithm that we use is not detailed here for the sake of
brevity. The reader is referred to Section 4
of~\cite{feliot2017bayesian} for a discussion about this
aspect. Details about the computation of the normalizing constant
$Z_n$ and about the variance of the proposed estimator are given in
Appendix~\ref{sec:covariance}.

\section{Numerical experiments}
\label{sec:experiments}

In our experiments, we illustrate the operation of the EWHI criterion
on the bi-objective BNH problem as defined in
\cite{chafekar2003constrained} for the following two weight functions
adapted from~\cite{zitzler2007hypervolume}:
\begin{equation}
  \left\{
    \begin{array}{lcl}
  \omega_1(y_1,y_2) &=& \displaystyle \frac{1}{15}e^{-\frac{y_{1}}{15}} \cdot \frac{\one_{[0,150]}(y_{1})}{150} \cdot \frac{\one_{[0,60]}(y_{2})}{60}\,,  \\[8pt]
      \omega_2(y_1,y_2) &=& \displaystyle \frac{1}{2}\left(\varphi\left(y,\mu_1,C\right) + \varphi\left(y,\mu_2,C\right)\right) \,,
    \end{array}
  \right.
\end{equation}
where $\varphi(y,\mu,C)$ denotes the Gaussian probability density
function with mean $\mu$ and covariance matrix $C$, evaluated at
$y$. The $\omega_1$ weight function is based on an exponential
distribution and encodes preference for the minimization of the first
objective. The $\omega_2$ weight function is a sum of two bivariate
Gaussian distributions and encodes preference for improving upon two
reference points $\mu_1$ and $\mu_2$, chosen as $\mu_1 = (80,20)$ and
$\mu_2 = (30,40)$ with $C = RS(RS)^T$, where
\begin{equation}
  \begin{array}{lcl}
    R = \left[
    \begin{array}{cc}
      \cos\left(\frac{\pi}{4}\right) & -\sin\left(\frac{\pi}{4}\right) \\
      \sin\left(\frac{\pi}{4}\right) & \cos\left(\frac{\pi}{4}\right)
    \end{array}
	\right]
&\mathrm{and}&
S = \left[
	\begin{array}{cc}
	20 & 0 \\
	0 & 3
	\end{array}
	\right]\,.
\end{array}
\end{equation}

To carry out the experiments, we use the BMOO algorithm
of~\cite{feliot2017thesis} with $m_x = m_y = 1000$ particles for both
SMC algorithms. The functions of the problem are modeled using
stationnary Gaussian processes with a constant mean and an anisotropic
Mat{\'e}rn covariance kernel. A log-normal prior distribution is
placed on the parameters of the kernel and these are updated at each
iteration of the algorithm using maximum a posteriori substition (see,
e.g., \cite{bect2012sequential}). The algorithm is initialized with a
pseudo-maximin latin hypercube design of $N=10$ experiments and is
iterated over 20 iterations. To handle the constraints of the BNH
problem, the EWHI criterion is multiplied by the probability of
feasibility, as is common practice in the Bayesian optimization
litterature (see, e.g., \cite{schonlau1998global}).

\begin{figure}
  \centering \setlength{\vSki}{5mm}
  \begin{minipage}[b]{5.5cm} \centering
  	\begin{psfrags}
  \psfrag{Objective space (N = 10)}[b][b]{} 
  \psfrag{f1}[t][t]{$f_1$} \psfrag{f2}[b][b]{\raisebox{2mm}{$f_2$}}
  \psfrag{0}[t][t]{\scriptsize 0} \psfrag{50}[t][t]{\scriptsize 50}
  \psfrag{100}[t][t]{\scriptsize 100} \psfrag{150}[t][t]{\scriptsize 150}
  \psfrag{0}[r][r]{\scriptsize 0} \psfrag{10}[r][r]{\scriptsize 10} \psfrag{20}[r][r]{\scriptsize 20} \psfrag{30}[r][r]{\scriptsize 30}
  \psfrag{40}[r][r]{\scriptsize 40} \psfrag{50}[r][r]{\scriptsize 50} \psfrag{60}[r][r]{\scriptsize 60}
    \includegraphics[width=5.5cm]{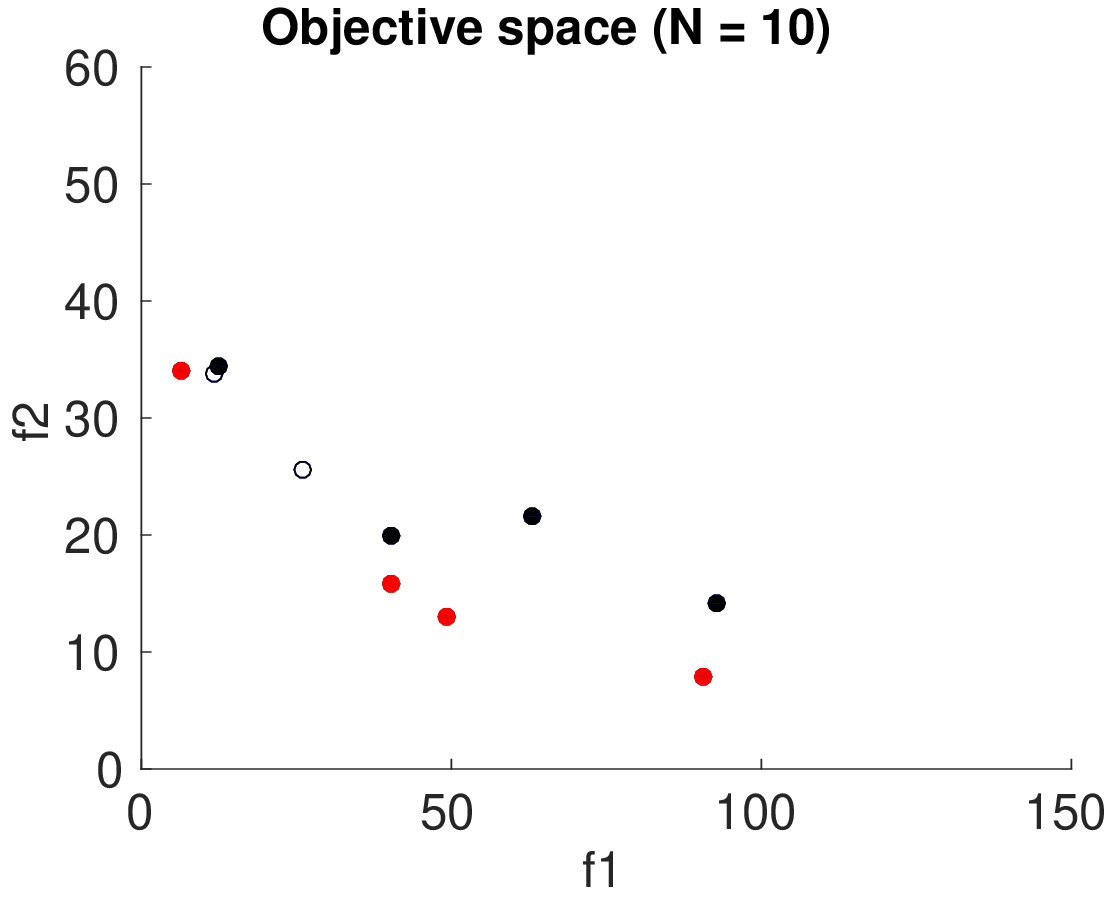} \\[1mm] 
    (a) EHVI ($N=10$)
  	\end{psfrags}
  \end{minipage} 
  \hspace{5mm}
  \begin{minipage}[b]{5.5cm} \centering
  	\begin{psfrags}
  \psfrag{Objective space (N = 30)}[b][b]{} 
  \psfrag{f1}[t][t]{$f_1$} \psfrag{f2}[b][b]{\raisebox{2mm}{$f_2$}}
  \psfrag{0}[t][t]{\scriptsize 0} \psfrag{50}[t][t]{\scriptsize 50}
  \psfrag{100}[t][t]{\scriptsize 100} \psfrag{150}[t][t]{\scriptsize 150}
  \psfrag{0}[r][r]{\scriptsize 0} \psfrag{10}[r][r]{\scriptsize 10} \psfrag{20}[r][r]{\scriptsize 20} \psfrag{30}[r][r]{\scriptsize 30}
  \psfrag{40}[r][r]{\scriptsize 40} \psfrag{50}[r][r]{\scriptsize 50} \psfrag{60}[r][r]{\scriptsize 60}
    \includegraphics[width=5.5cm]{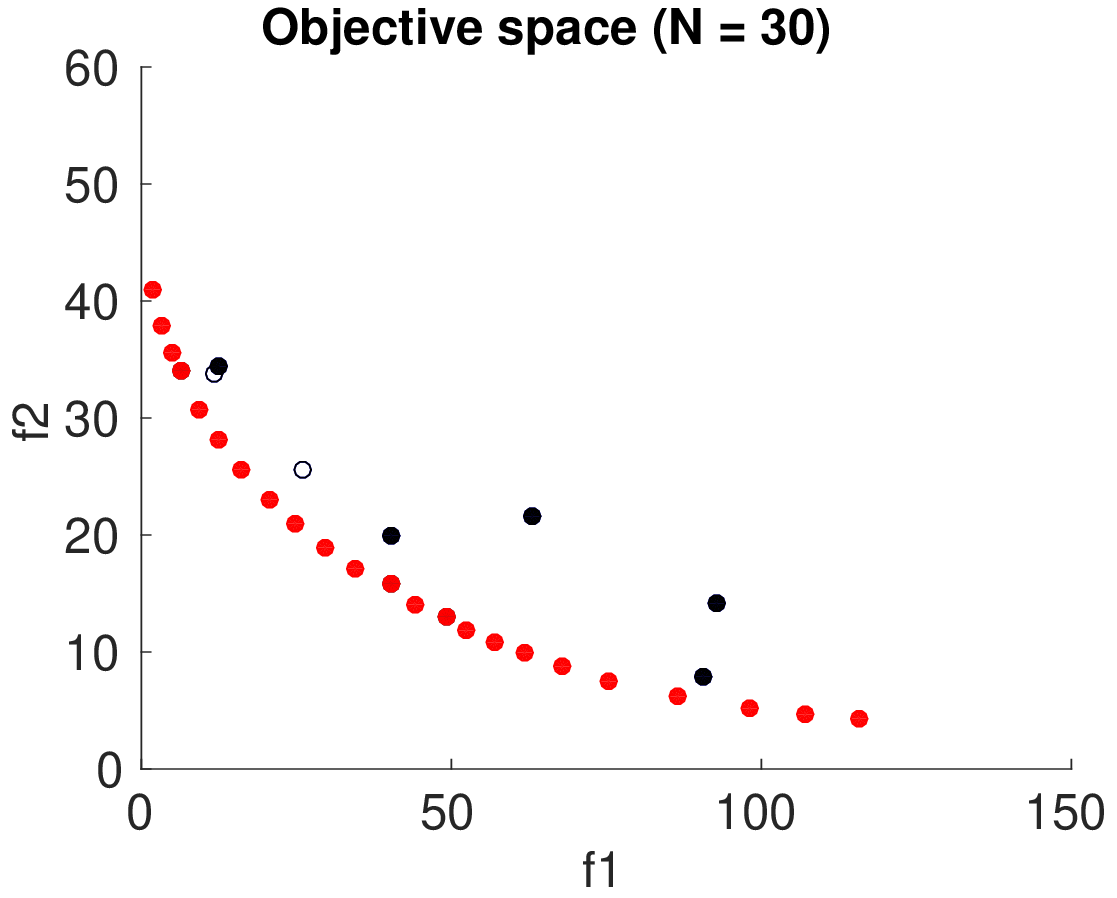} \\[1mm]
    (b) EHVI ($N=30$)
  	\end{psfrags}
  \end{minipage}\\[\vSki]
  \begin{minipage}[b]{5.5cm} \centering
  	\begin{psfrags}
  \psfrag{Objective space (N = 10)}[b][b]{} 
  \psfrag{f1}[t][t]{$f_1$} \psfrag{f2}[b][b]{\raisebox{2mm}{$f_2$}}
  \psfrag{0}[t][t]{\scriptsize 0} \psfrag{50}[t][t]{\scriptsize 50}
  \psfrag{100}[t][t]{\scriptsize 100} \psfrag{150}[t][t]{\scriptsize 150}
  \psfrag{0}[r][r]{\scriptsize 0} \psfrag{10}[r][r]{\scriptsize 10} \psfrag{20}[r][r]{\scriptsize 20} \psfrag{30}[r][r]{\scriptsize 30}
  \psfrag{40}[r][r]{\scriptsize 40} \psfrag{50}[r][r]{\scriptsize 50} \psfrag{60}[r][r]{\scriptsize 60}
    \includegraphics[width=5.5cm]{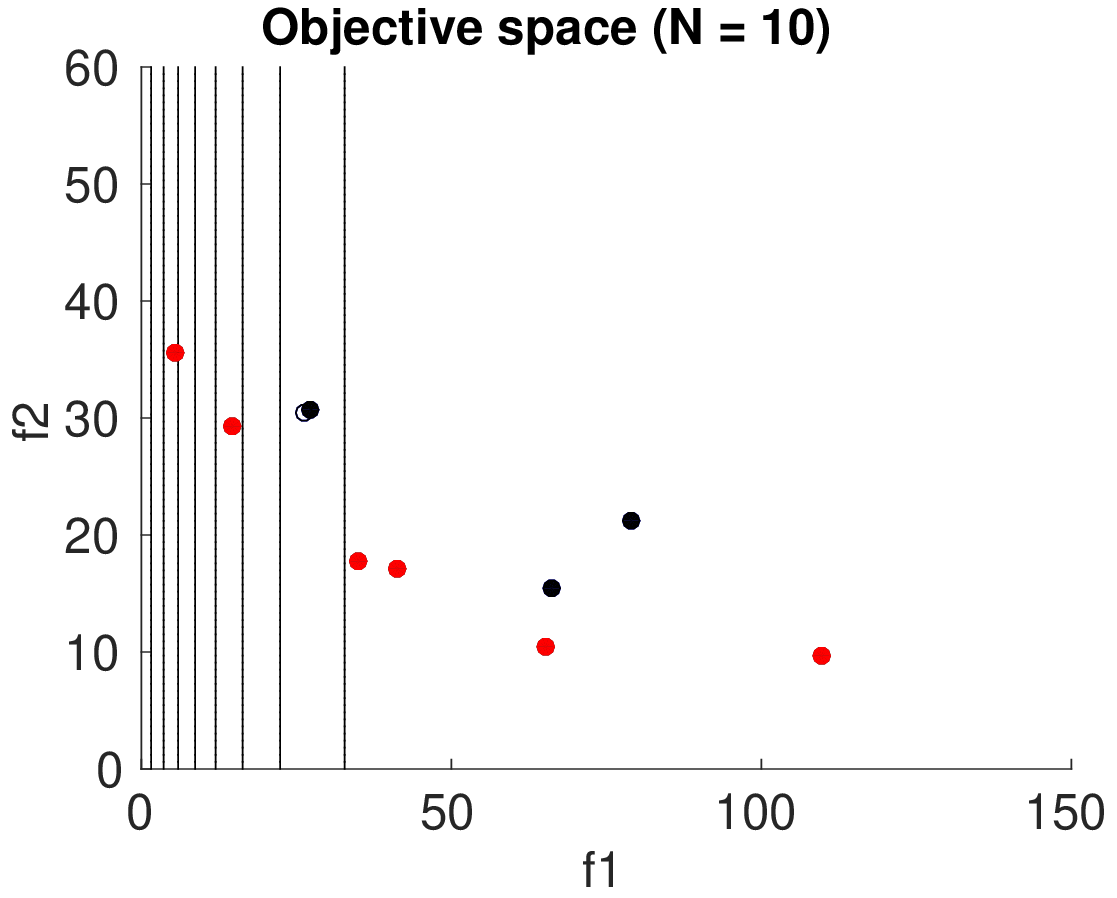} \\[1mm] 
    (c) EWHI with $\omega_1$ ($N=10$)
  	\end{psfrags}
  \end{minipage} 
  \hspace{5mm}
  \begin{minipage}[b]{5.5cm} \centering
  	\begin{psfrags}
  \psfrag{Objective space (N = 30)}[b][b]{} 
  \psfrag{f1}[t][t]{$f_1$} \psfrag{f2}[b][b]{\raisebox{2mm}{$f_2$}}
  \psfrag{0}[t][t]{\scriptsize 0} \psfrag{50}[t][t]{\scriptsize 50}
  \psfrag{100}[t][t]{\scriptsize 100} \psfrag{150}[t][t]{\scriptsize 150}
  \psfrag{0}[r][r]{\scriptsize 0} \psfrag{10}[r][r]{\scriptsize 10} \psfrag{20}[r][r]{\scriptsize 20} \psfrag{30}[r][r]{\scriptsize 30}
  \psfrag{40}[r][r]{\scriptsize 40} \psfrag{50}[r][r]{\scriptsize 50} \psfrag{60}[r][r]{\scriptsize 60}
    \includegraphics[width=5.5cm]{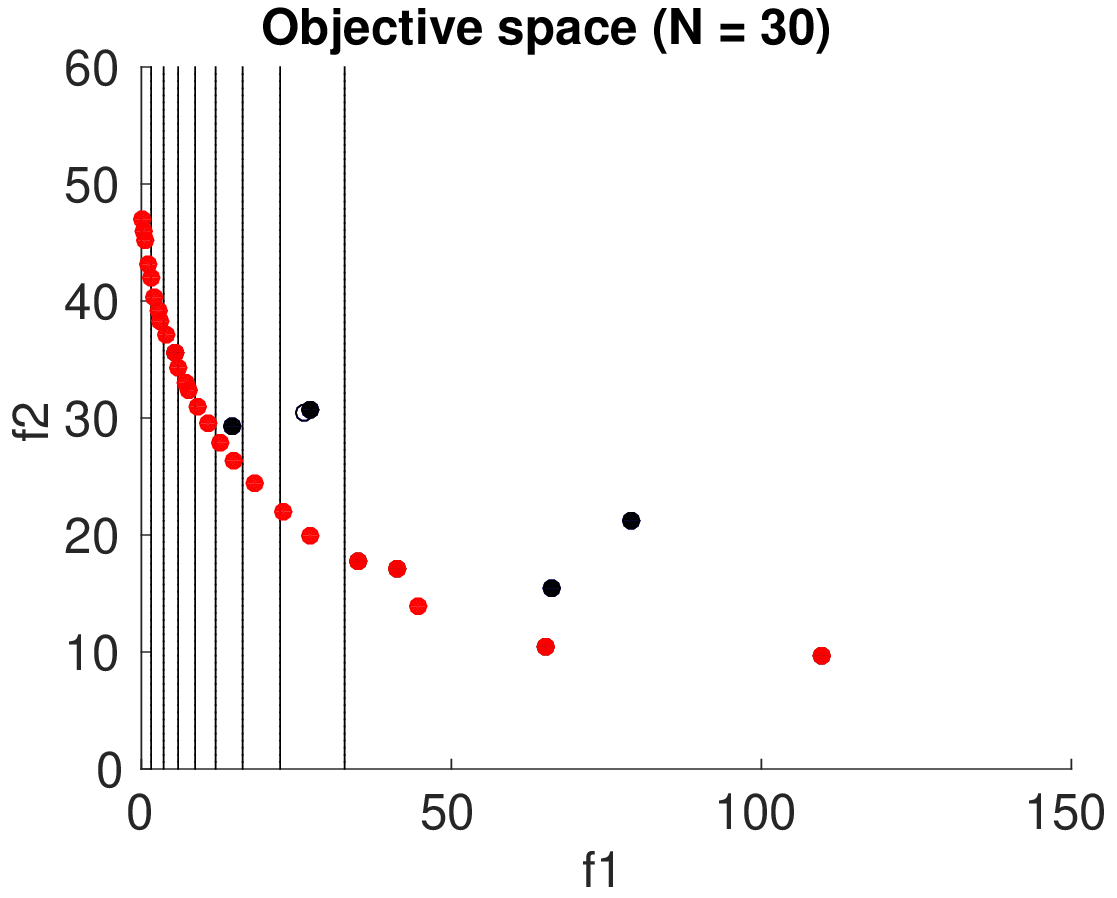} \\[1mm]
    (d) EWHI with $\omega_1$ ($N=30$)
  	\end{psfrags}
  \end{minipage}\\[\vSki]
  \begin{minipage}[b]{5.5cm} \centering
  	\begin{psfrags}
  \psfrag{Objective space (N = 10)}[b][b]{} 
  \psfrag{f1}[t][t]{$f_1$} \psfrag{f2}[b][b]{\raisebox{2mm}{$f_2$}}
  \psfrag{0}[t][t]{\scriptsize 0} \psfrag{50}[t][t]{\scriptsize 50}
  \psfrag{100}[t][t]{\scriptsize 100} \psfrag{150}[t][t]{\scriptsize 150}
  \psfrag{0}[r][r]{\scriptsize 0} \psfrag{10}[r][r]{\scriptsize 10} \psfrag{20}[r][r]{\scriptsize 20} \psfrag{30}[r][r]{\scriptsize 30}
  \psfrag{40}[r][r]{\scriptsize 40} \psfrag{50}[r][r]{\scriptsize 50} \psfrag{60}[r][r]{\scriptsize 60}
    \includegraphics[width=5.5cm]{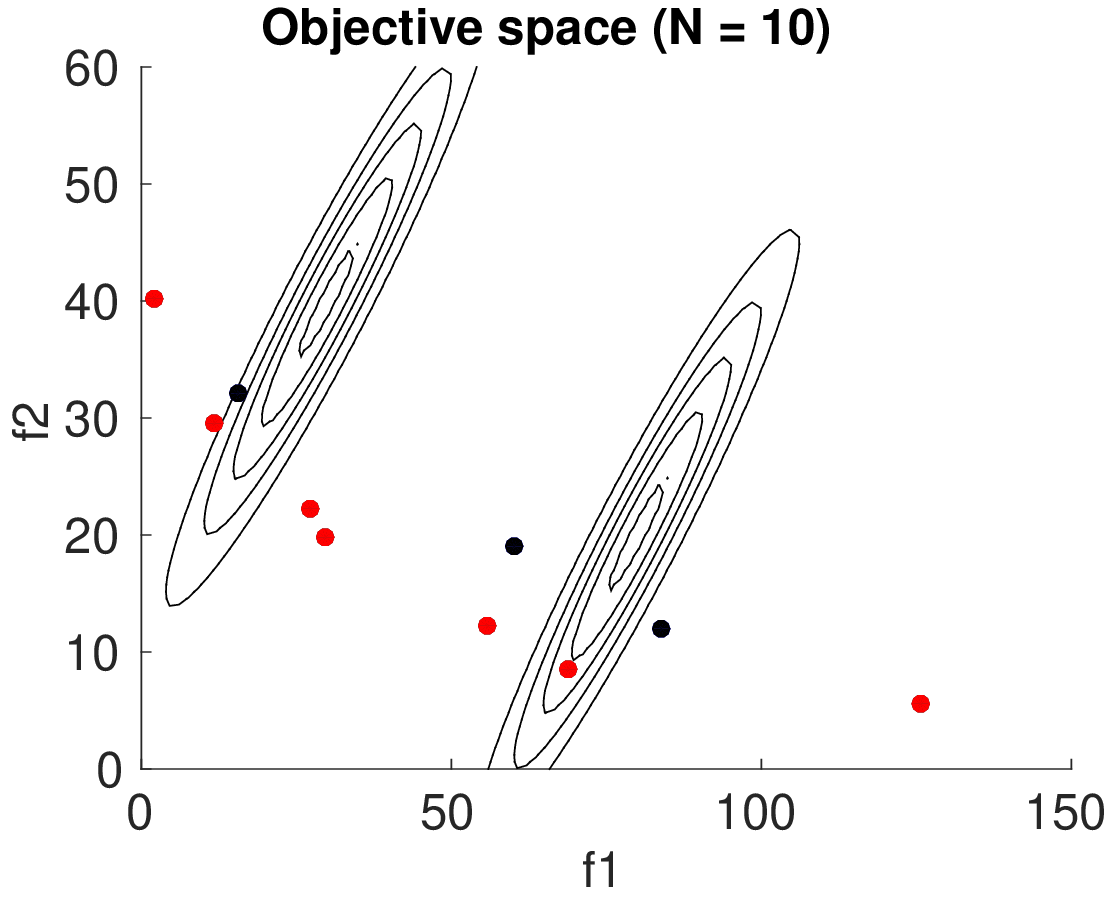} \\[1mm] 
    (e) EWHI with $\omega_2$ ($N=10$)
  	\end{psfrags}
  \end{minipage} 
  \hspace{5mm}
  \begin{minipage}[b]{5.5cm} \centering
  	\begin{psfrags}
  \psfrag{Objective space (N = 30)}[b][b]{} 
  \psfrag{f1}[t][t]{$f_1$} \psfrag{f2}[b][b]{\raisebox{2mm}{$f_2$}}
  \psfrag{0}[t][t]{\scriptsize 0} \psfrag{50}[t][t]{\scriptsize 50}
  \psfrag{100}[t][t]{\scriptsize 100} \psfrag{150}[t][t]{\scriptsize 150}
  \psfrag{0}[r][r]{\scriptsize 0} \psfrag{10}[r][r]{\scriptsize 10} \psfrag{20}[r][r]{\scriptsize 20} \psfrag{30}[r][r]{\scriptsize 30}
  \psfrag{40}[r][r]{\scriptsize 40} \psfrag{50}[r][r]{\scriptsize 50} \psfrag{60}[r][r]{\scriptsize 60}
    \includegraphics[width=5.5cm]{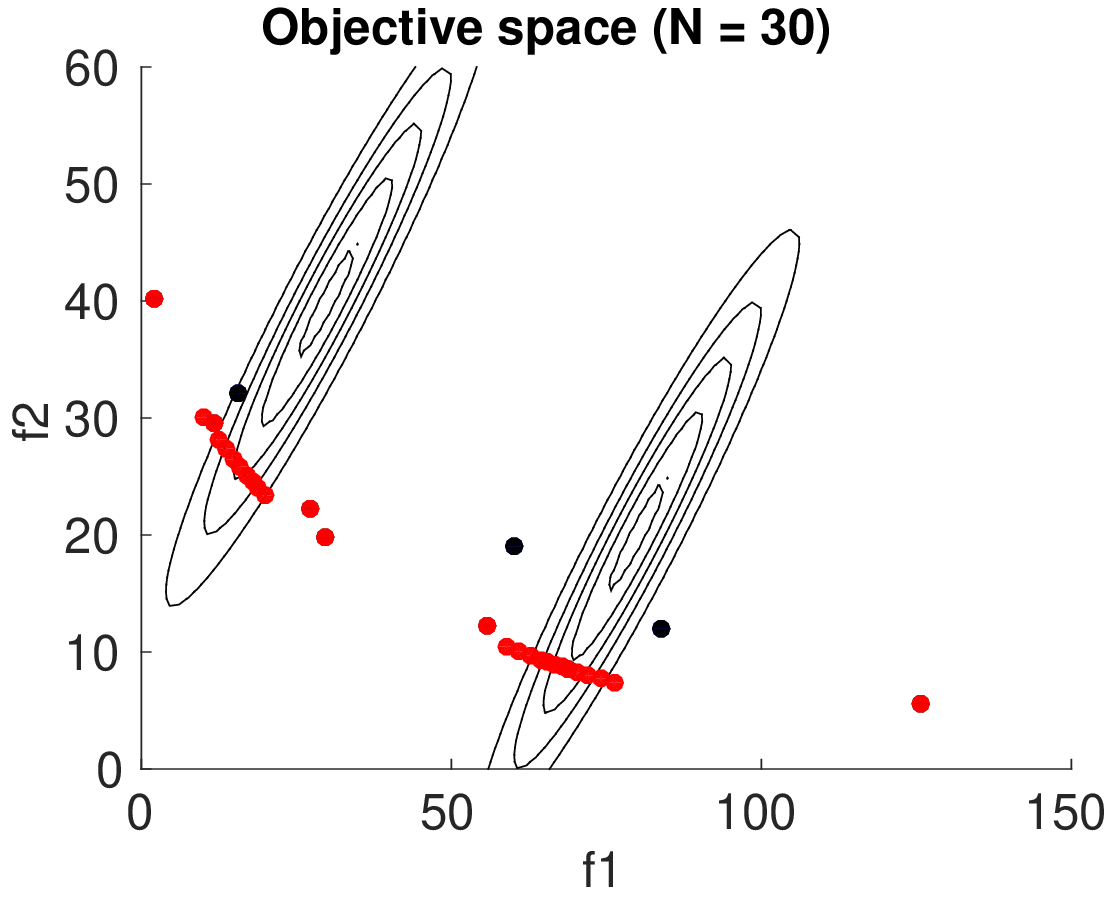} \\[1mm]
    (f) EWHI with $\omega_2$ ($N=30$)
  	\end{psfrags}
  \end{minipage}
  \caption{Distributions obtained after 20 iterations of the
    optimization algorithm on the BNH problem when the weight
    functions $\omega_1$ and $\omega_2$ are used. The results obtained
    using the EHVI criterion are shown for reference. The contours of
    the weight functions are represented as black lines and the
    non-dominated solutions as red disks. Black disks indicate
    feasible dominated solutions and black circles indicate
    non-feasible solutions.}
  \label{fig:illu_preferences}
\end{figure}

In Figure~\ref{fig:illu_preferences}, results obtained by the
algorithm using the weight functions $\omega_1$ and $\omega_2$ in the
EWHI definition are compared to results obtained by the same algorithm
using the EHVI criterion. Observe in
Figures~\ref{fig:illu_preferences}(d)
and~\ref{fig:illu_preferences}(f) that observations are concentrated
in regions of the Pareto front that correspond to high $\omega$
values, whereas observations are spread along the front in
Figure~\ref{fig:illu_preferences}(b) where the EHVI is used. In
practice, this means that less iterations would have been required to
satisfyingly populate the interesting regions of the Pareto front.

\section{Conclusions and perspectives}
\label{sec:conclusions}

It is shown in this paper how user-defined weight functions can be
leveraged by a Bayesian framework to produce optimization strategies
that focus on preferred regions of the Pareto front of multi-objective
optimization problems. Two example weight functions
from~\cite{zitzler2007hypervolume} which encode respectively a
preference for one objective and a preference toward specific regions
of the Pareto front are used, and the demonstration of the
effectiveness of the proposed approach is carried out on a simple
bi-objective optimization problem.

On more practical problems, crafting sensible weight functions can be
a difficult task, especially when one has no prior knowledge about the
approximate location of the Pareto front. The use of desirability
functions (see,
e.g. \cite{emmerich2014reference,harrington1965desirability,wagner2010integration})
or utility functions (see, e.g., \cite{astudillomulti2017}) might
provide useful insights on that issue and shall be the object of
future investigations to provide a more principled approach.

In the presented framework, optimization strategies are built
sequentially using an expected improvement sampling criterion called
the expected weighted hypervolume improvement (EWHI) criterion. The
exact computation of the criterion being intractable in general, an
approximate computation prodecure using importance sampling is
proposed. A sampling density that is optimal for the simultaneous
computation of the criterion for a set of candidate points is crafted
and a sequential Monte-Carlo algorithm is used to produce samples from
this density.

This choice triggers an immediate question: What is the sample size
$m_y$ required by the algorithm? In fact, the problem is not so much
to obtain a precise approximation of $\EI_n$ for all $x\in \XX_n$,
which would require a large sample size to distinguish very close
points, but to deal with the optimization
problem~\eqref{eq:oneStepLookAhead} and to identify with good
confidence the points of $\XX_n$ that correspond to high values of
$\EI_n$. A first step toward a solution to this problem is to compute
an approximation of the variance of $\hat\EI_n$, as carried out in
Appendix~\ref{sec:covariance}. Further investigations on this issue
are left for future work.

\appendix
\section{Approximate variance of the EI estimator}
\label{sec:covariance}

We derive in this appendix the variance of the SMC estimator for
$\rho_n$. In the SMC procedure that we consider, the particles
$\left(y_{n,i}\right)_{1 \leq i \leq m}$ are obtained from a sequence
of densities $(\pi_{n,t})_{0 \leq t \leq T}$, where $\pi_{n,0}$ is an
easy-to-sample initial density and $\pi_{n,T} = \pi_n$ is the target
density. Let $(\gamma_{n,t})_{0 \leq t \leq T}$ and
$(Z_{n,t})_{0 \leq t \leq T}$ denote the corresponding sequences of
un-normalized densities and normalizing constants.

First, observe that, for $1\leq t \leq T$,
\begin{equation}
\label{eq:gamma-recur}
\begin{array}{lcl}
Z_{n,t} &=& \displaystyle \int_{\Hnc} \gamma_{n,t}(y)\,\dy\\
        &=& \displaystyle Z_{n,t-1} \int_{G_n} \frac{\gamma_{n,t}(y)}{\gamma_{n,t-1}(y)}\, \pi_{n,t-1}(y)\, \dy \,.
\end{array}
\end{equation}
Thus, we can derive a sequence of approximations $\hat{Z}_{n,t}$ of
$Z_{n,t}$, $t\geq1$, using the following recursion formula:
\begin{equation}
  \label{eq:Gamma-approx}
  \left\{
    \begin{array}{l}
      \hat Z_{n,0} =  Z_{n,0} = \int_{G_n} \gamma_{n,0}(y)\, \dy,\\
      \hat{Z}_{n,t} = \hat{Z}_{n,t-1} \left( \frac{1}{m} \sum_{i=1}^{m} \frac{\gamma_{n,t}(y_{n,t-1,i})}{\gamma_{n,t-1}(y_{n,t-1,i})} \right) \,,  
    \end{array}\right.
\end{equation}
where the particles $(y_{n,t-1,i})_{1 \leq i \leq m}\sim \pi_{n,t-1}$
are obtained using an SMC procedure (see, e.g.,
\cite{bect2017bayesian}). The estimator of~$\rho_n(x)$ that we
actually consider is then
\begin{equation}
  \label{eq:approx-ei-product}
  \hat{\rho}_n(x) = \displaystyle \frac{\hat{Z}_n}{m} \sum_{i=1}^m
  \frac{\omega(y)\,\P_n \left( \xi(x) \prec y_{n,i} \right)}{\gamma_n(y_{n,i})} =
  \hat Z_n \hat\alpha_{n}(x)
\end{equation}
where 
\begin{equation}
  \hat{\alpha}_n(x) = \frac{1}{m} \sum_{i=1}^m \frac{\omega(y)\,\P_n \left( \xi(x) \prec y_{n,i} \right)}{\gamma_n(y_{n,i})}\,,
\end{equation}
and
\begin{equation}
  \hat Z_n = \hat Z_{n,T} = Z_{n,0} \prod_{u=1}^T \hat{\theta}_{n,\,u}\,,
\end{equation}
with
\begin{equation}
  \hat{\theta}_{n,t} = \frac{1}{m} \sum_{i=1}^{m} \frac{\gamma_{n,t}(y_{n,t-1,i})}{\gamma_{n,t-1}(y_{n,t-1,i})}\,.
\end{equation}

Now, assume the idealized setting, as usual in the SMC literature
(see, e.g., \cite{cerou2012sequential}), where
\begin{enumerate}[(i)]
\item $y_{n,t,i} \stackrel{\scriptstyle
    \text{i.i.d}}{\sim} \pi_{n,t}$, $1 \leq i \leq m$,
\item the samples $\YY_{n,t} = (y_{n,t,i})_{1 \leq i \leq m}$ are
  independent, $0\leq t \leq T$.
\end{enumerate}

Observe from~\eqref{eq:ei_int_approx} and~\eqref{eq:gamma-recur} that
under~(i), $\hat{\alpha}_n(x)$ is an unbiased estimator of
$\alpha_n(x) = \frac{\rho_n(x)}{Z_n}$, and $\hat{\theta}_{n,t}$ is an
unbiased estimator of $\theta_{n,t} = \frac{Z_{n,t}}{Z_{n,t-1}}$,
$1 \leq t \leq T$. Moreover, under~(ii), $\hat{\alpha}_n(x)$ and the
$(\hat{\theta}_{n,t})_{1 \leq t \leq T}$ are independent. Thus,
\begin{equation*}
\begin{array}{l}
  \mathrm{Var}\, \hat{\rho}_n(x) 
  = \E \bigl( \hat{\alpha}_n^2 \bigr) \E \bigl( \hat Z_n^2 \bigr) - \E \bigl( \hat{\alpha}_n(x) \bigr)^2  \E \bigl( \hat Z_n \bigr)^2 \\[4pt]
  = \bigl( \mathrm{Var}\,  \hat{\alpha}_n(x)\!  +\! \alpha_n(x)^2 \bigr)  \bigl( \mathrm{Var}\,  \hat Z_n  + Z_n^2 \bigr)  - \alpha_n(x)^2  Z_n^2 \\[4pt]
  = \mathrm{Var}\, \hat{\alpha}_n(x)  \mathrm{Var}\,  \hat Z_n + \alpha_n(x)^2  \mathrm{Var}\,  \hat Z_n + Z_n^2  \mathrm{Var} \, \hat\alpha_n(x) \,. 
\end{array}
\end{equation*}
We obtain the coefficient of variation of $\hat\EI_n(x)$
\begin{equation}
\label{eq:lambda-recur}
\frac{\mathrm{Var}\,  \hat\rho_n(x) }{\rho_n(x)^2} = \Lambda_n(x)^2 + \left( 1+\Lambda_n(x)^2 \right)  \Delta_{n,T}^2,
\end{equation}
where
$\Lambda_n(x)^{2} =
\frac{\mathrm{Var}\,\hat\alpha_n(x)}{\alpha_n(x)^2}$ and
$\Delta_{n,t}^2 = \frac{\mathrm{Var} \hat Z_{n,t}}{Z_{n,t}^2}$ are the
coefficients of variation of $\hat\alpha_n(x)$ and $\hat Z_{n,t}$
respectively.

Using the same ideas as above, we have 
\begin{equation}
\label{eq:delta-recur}
\Delta_{n,t}^2 = \delta_{n,t}^2 + \left( 1+\delta_{n,t}^2 \right) \Delta_{n,t-1}^2,
\end{equation}
where
$\delta_{n,t}^2 = \frac{\mathrm{Var}\,
  \hat{\theta}_{n,t}}{\theta_{n,t}^2}$ is the coefficient of variation
of $\hat{\theta}_{n,t}$.

Estimators of $\Lambda_{n}(x)^2$, $\Delta_{n,t}^2$ and
$\delta_{n,t}^2$ can be derived under~(ii). For instance, observe that
\begin{equation}
  \delta_{n,t}^2 = \displaystyle \frac{1}{m}\,\frac{\mathrm{Var}\biggl( \frac{\gamma_{n,t}(y_{n,t-1,\,1})}{\gamma_{n,t-1}(y_{n,t-1,\,1})} \biggr)}{\E\biggl( \frac{\gamma_{n,t}(y_{n,t-1,\,1})}{\gamma_{n,t-1}(y_{n,t-1,\,1})} \biggr)^2}\,.
\end{equation}
Thus, an estimator of $\delta_{n,t}^2$ is
\begin{equation}
  \label{eq:deltant}
  \hat{\delta}_{n,t}^2 = \displaystyle \frac{ \sum_{i=1}^{m} \frac{\gamma_{n,t}(y_{n,t-1,i})^2}{\gamma_{n,t-1}(y_{n,t-1,i})^2} }{\biggl( \sum_{i=1}^{m} \frac{\gamma_{n,t}(y_{n,t-1,i})}{\gamma_{n,t-1}(y_{n,t-1,i})} \biggr)^2} - \frac{1}{m}.
\end{equation}

Plugging (\ref{eq:deltant}) in~(\ref{eq:delta-recur}), we obtain an
estimator of $\Delta_{n,t}^2$:
\begin{equation}
\label{eq:Delta-approx}
\hat{\Delta}_{n,t}^2 = \hat{\delta}_{n,t}^2 + \left( 1+\hat{\delta}_{n,t}^2 \right) \cdot \hat{\Delta}_{n,t-1}^2.
\end{equation}

Similarly, an estimator of $\Lambda_n(x)^2$ is
\begin{equation}
\label{eq:Lambda-approx}
\hat\Lambda_n(x)^2 = \displaystyle \frac{ \sum_{i=1}^{m} \frac{\omega(y)^2\, \P_n \left( \xi(x) \prec y_{n,i} \right)^2}{\gamma_n(y_{n,i})^2}}{\left( \sum_{i=1}^{m} \frac{\omega(y)\, \P_n \left( \xi(x) \prec y_{n,i} \right)}{\gamma_n(y_{n,i})} \right)^2} - \frac{1}{m}.
\end{equation}

As a result, we obtain the following numerically tractable
approximation of the variance of $\hat\rho_n(x)$:
\begin{equation}
\label{eq:var-approx}
\mathrm{Var} \left( \hat\rho_n(x)  \right)
\approx \displaystyle \hat\rho_n(x)^2 \cdot \left( \hat{\Lambda}_n(x)^2 + \left( 1+\hat{\Lambda}_n(x)^2 \right) \cdot \hat{\Delta}_{n,t}^2 \right),
\end{equation}
where $\hat{Z}_{n,t}$ and $\hat{\Delta}_{n,t}^2$ are obtained
recursively using~(\ref{eq:Gamma-approx}) and~(\ref{eq:Delta-approx}),
$\hat{\Lambda}_n(x)^2$ is computed using~(\ref{eq:Lambda-approx}) and
$\hat\rho_n(x)$ is computed using~(\ref{eq:approx-ei-product}).

\bibliographystyle{splncs04}
\bibliography{LOD2018-EWHI-paper}

\begin{thebibliography}{10}
\providecommand{\url}[1]{\texttt{#1}}
\providecommand{\urlprefix}{URL }
\providecommand{\doi}[1]{https://doi.org/#1}

\bibitem{astudillomulti2017}
Astudillo, R., Frazier, P.: Multi-attribute bayesian optimization under utility
  uncertainty. Proceedings of the NIPS Workshop on Bayesian Optimization,
  December 2017, Long Beach, USA  (To appear)

\bibitem{au2001estimation}
Au, S.K., Beck, J.L.: Estimation of small failure probabilities in high
  dimensions by subset simulation. Probabilistic Engineering Mechanics
  \textbf{16}(4),  263--277 (2001)

\bibitem{auger2009articulating}
Auger, A., Bader, J., Brockhoff, D., Zitzler, E.: Articulating user preferences
  in many-objective problems by sampling the weighted hypervolume. In:
  Proceedings of the 11th Annual conference on Genetic and evolutionary
  computation. pp. 555--562. ACM (2009)

\bibitem{auger2009investigating}
Auger, A., Bader, J., Brockhoff, D., Zitzler, E.: Investigating and exploiting
  the bias of the weighted hypervolume to articulate user preferences. In:
  Proceedings of the 11th Annual conference on Genetic and evolutionary
  computation. pp. 563--570. ACM (2009)

\bibitem{bect2012sequential}
Bect, J., Ginsbourger, D., Li, L., Picheny, V., Vazquez, E.: Sequential design
  of computer experiments for the estimation of a probability of failure.
  Statistics and Computing  \textbf{22}(3),  773--793 (2012)

\bibitem{bect2017bayesian}
Bect, J., Li, L., Vazquez, E.: Bayesian subset simulation. SIAM/ASA Journal on
  Uncertainty Quantification  \textbf{5}(1),  762--786 (2017)

\bibitem{bect2015echantillonnage}
Bect, J., Sueur, R., G{\'e}rossier, A., Mongellaz, L., Petit, S., Vazquez, E.:
  {\'E}chantillonnage pr{\'e}f{\'e}rentiel et m{\'e}ta-mod{\`e}les:
  m{\'e}thodes bay{\'e}siennes optimale et d{\'e}fensive. In: 47{\`e}mes
  Journ{\'e}es de Statistique de la SFdS-JdS 2015 (2015)

\bibitem{benassi2013nouvel}
Benassi, R.: Nouvel algorithme d'optimisation bay{\'e}sien utilisant une
  approche Monte-Carlo s{\'e}quentielle. Ph.D. thesis, Sup{\'e}lec (2013)

\bibitem{cerou2012sequential}
C{\'e}rou, F., Del~Moral, P., Furon, T., Guyader, A.: Sequential {M}onte
  {C}arlo for rare event estimation. Statistics and Computing  \textbf{22}(3),
  795--808 (2012)

\bibitem{chafekar2003constrained}
Chafekar, D., Xuan, J., Rasheed, K.: Constrained multi-objective optimization
  using steady state genetic algorithms. In: Genetic and Evolutionary
  Computation-GECCO 2003. pp. 813--824. Springer (2003)

\bibitem{del2006sequential}
Del~Moral, P., Doucet, A., Jasra, A.: Sequential monte carlo samplers. Journal
  of the Royal Statistical Society: Series B (Statistical Methodology)
  \textbf{68}(3),  411--436 (2006)

\bibitem{emmerich2005}
Emmerich, M.: Single- and multiobjective evolutionary design optimization
  assisted by {G}aussian random field metamodels. Ph.D. thesis, Technical
  University Dortmund (2005)

\bibitem{emmerich2014reference}
Emmerich, M., Deutz, A.H., Yevseyeva, I.: On reference point free weighted
  hypervolume indicators based on desirability functions and their
  probabilistic interpretation. Procedia Technology  \textbf{16},  532--541
  (2014)

\bibitem{emmerich2006single}
Emmerich, M., Giannakoglou, K.C., Naujoks, B.: Single- and multi-objective
  evolutionary optimization assisted by {Ga}ussian random field metamodels.
  IEEE Transactions on Evolutionary Computation  \textbf{10}(4),  421--439
  (2006)

\bibitem{emmerich2008computation}
Emmerich, M., Klinkenberg, J.W.: The computation of the expected improvement in
  dominated hypervolume of {P}areto front approximations. Technical report,
  Leiden University  (2008)

\bibitem{feliot2017thesis}
Feliot, P.: A Bayesian approach to constrained multi-objective optimization.
  Ph.D. thesis, IRT SystemX and Centrale Sup{\'e}lec (2017)

\bibitem{feliot2017bayesian}
Feliot, P., Bect, J., Vazquez, E.: A {B}ayesian approach to constrained
  single-and multi-objective optimization. Journal of Global Optimization
  \textbf{67}(1-2),  97--133 (2017)

\bibitem{GinsbLeRiche2009}
Ginsbourger, D., Le~Riche, R.: Towards {G}aussian process-based optimization
  with finite time horizon. In: mODa 9 -- Advances in Model-Oriented Design and
  Analysis, pp. 89--96. Springer (2010)

\bibitem{harrington1965desirability}
Harrington, E.C.: The desirability function. Industrial quality control
  \textbf{21}(10),  494--498 (1965)

\bibitem{jones1998efficient}
Jones, D.R., Schonlau, M., Welch, W.J.: Efficient global optimization of
  expensive black-box functions. Journal of Global Optimization
  \textbf{13}(4),  455--492 (1998)

\bibitem{knowles2002metrics}
Knowles, J., Corne, D.: On metrics for comparing nondominated sets. In:
  Proceedings of the 2002 Congress on Evolutionary Computation, 2002. CEC'02.
  vol.~1, pp. 711--716. IEEE (2002)

\bibitem{kushner1964new}
Kushner, H.J.: A new method of locating the maximum point of an arbitrary
  multipeak curve in the presence of noise. Journal of Fluids Engineering
  \textbf{86}(1),  97--106 (1964)

\bibitem{laumanns1999approximating}
Laumanns, M., Rudolph, G., Schwefel, H.P.: Approximating the pareto set:
  Concepts, diversity issues, and performance assessment. Secretary of the SFB
  531 (1999)

\bibitem{mockus78}
Mockus, J., Tiesis, V., \v{Z}ilinskas, A.: The application of {B}ayesian
  methods for seeking the extremum. In: Dixon, L.C.W., Szeg{\"{o}}, G.P. (eds.)
  Towards Global Optimization. vol.~2, pp. 117--129. North Holland, New York
  (1978)

\bibitem{santner2003design}
Santner, T.J., Williams, B.J., Notz, W.: The design and analysis of computer
  experiments. Springer (2003)

\bibitem{schonlau1998global}
Schonlau, M., Welch, W.J., Jones, D.R.: Global versus local search in
  constrained optimization of computer models. In: New Developments and
  Applications in Experimental Design: Selected Proceedings of a 1997 Joint
  AMS-IMS-SIAM Summer Conference. IMS Lecture Notes-Monographs Series, vol.~34,
  pp. 11--25. Institute of Mathematical Statistics (1998)

\bibitem{wagner2010integration}
Wagner, T., Trautmann, H.: Integration of preferences in hypervolume-based
  multiobjective evolutionary algorithms by means of desirability functions.
  IEEE Transactions on Evolutionary Computation  \textbf{14}(5),  688--701
  (2010)

\bibitem{williams2006gaussian}
Williams, C.K.I., Rasmussen, C.: {G}aussian processes for machine learning. the
  MIT Press  \textbf{2}(3), ~4 (2006)

\bibitem{zitzler2007hypervolume}
Zitzler, E., Brockhoff, D., Thiele, L.: The hypervolume indicator revisited: On
  the design of pareto-compliant indicators via weighted integration. In:
  Evolutionary multi-criterion optimization. pp. 862--876. Springer (2007)

\bibitem{zitzler1998multiobjective}
Zitzler, E., Thiele, L.: Multiobjective optimization using evolutionary
  algorithms--a comparative case study. In: Parallel problem solving from
  nature--PPSN V. pp. 292--301. Springer (1998)

\end{thebibliography}
\end{document}